\magnification\magstep1
\input epsf

\font\tengoth=eufm10  \font\ninegoth=eufm9
\font\eightgoth=eufm8 \font\sevengoth=eufm7
\font\sixgoth=eufm6   \font\fivegoth=eufm5
\newfam\gothfam \def\goth{\fam\gothfam\tengoth}
\textfont\gothfam=\tengoth
\scriptfont\gothfam=\sevengoth
\scriptscriptfont\gothfam=\fivegoth

\font\ninerm=cmr9  \font\eightrm=cmr8  \font\sixrm=cmr6
\font\ninei=cmmi9  \font\eighti=cmmi8  \font\sixi=cmmi6
\font\ninesy=cmsy9 \font\eightsy=cmsy8 \font\sixsy=cmsy6
\font\ninebf=cmbx9 \font\eightbf=cmbx8 \font\sixbf=cmbx6
\font\nineit=cmti9 \font\eightit=cmti8
\font\ninett=cmtt9 \font\eighttt=cmtt8
\font\ninesl=cmsl9 \font\eightsl=cmsl8
\newskip\ttglue

% switch to 8-point type
\def\eightpoint{\def\rm{\fam0\eightrm}
  \textfont0=\eightrm \scriptfont0=\sixrm
  \scriptscriptfont0\fiverm
  \textfont1=\eighti \scriptfont1=\sixi
  \scriptscriptfont1\fivei
  \textfont2=\eightsy \scriptfont2=\sixsy
  \scriptscriptfont2\fivesy
  \textfont3=\tenex \scriptfont3=\tenex
  \scriptscriptfont3\tenex
  \textfont\itfam=\eightit\def\it{\fam\itfam\eightit}%
  \textfont\slfam=\eightsl\def\sl{\fam\slfam\eightsl}%
  \textfont\ttfam=\eighttt\def\tt{\fam\ttfam\eighttt}%
  \textfont\gothfam=\eightgoth\scriptfont\gothfam=\sixgoth
  \scriptscriptfont\gothfam=\fivegoth
  \def\goth{\fam\gothfam\tengoth}
  \textfont\bffam=\eightbf\scriptfont\bffam=\sixbf
  \scriptscriptfont\bffam=\fivebf
  \def\bf{\fam\bffam\eightbf}%
  \tt\ttglue=.5em plus.25em minus.15em
  \normalbaselineskip=9pt \setbox\strutbox\hbox{\vrule
  height7pt depth2pt width0pt}%
  \let\big=\eightbig\normalbaselines\rm}

% switch to 9-point type
\def\ninepoint{\def\rm{\fam0\ninerm}
  \textfont0=\ninerm \scriptfont0=\sixrm
  \scriptscriptfont0\fiverm
  \textfont1=\ninei \scriptfont1=\sixi
  \scriptscriptfont1\fivei
  \textfont2=\ninesy \scriptfont2=\sixsy
  \scriptscriptfont2\fivesy
  \textfont3=\tenex \scriptfont3=\tenex
  \scriptscriptfont3\tenex
  \textfont\itfam=\nineit\def\it{\fam\itfam\nineit}%
  \textfont\slfam=\ninesl\def\sl{\fam\slfam\ninesl}%
  \textfont\ttfam=\ninett\def\tt{\fam\ttfam\ninett}%
  \textfont\gothfam=\ninegoth\scriptfont\gothfam=\sixgoth
  \scriptscriptfont\gothfam=\fivegoth
  \def\goth{\fam\gothfam\tengoth}
  \textfont\bffam=\ninebf\scriptfont\bffam=\sixbf
  \scriptscriptfont\bffam=\fivebf
  \def\bf{\fam\bffam\ninebf}%
  \tt\ttglue=.5em plus.25em minus.15em
  \normalbaselineskip=11pt \setbox\strutbox\hbox{\vrule
  height8pt depth3pt width0pt}%
  \let\big=\ninebig\normalbaselines\rm}
\def\bibliography#1\par{\vskip0pt
  plus.3\vsize\penalty-250\vskip0pt
  plus-.3\vsize\bigskip\vskip\parskip
  \message{Bibliography}\leftline{\bf
  Bibliography}\nobreak\smallskip\noindent
  \ninepoint\frenchspacing#1}

\hyphenation{Prin-ce-ton}
\def\quid{$\!\!\!\!\!\quad$}
\def\ZZ{{\bf Z}}
\def\FF{{\bf F}}
\def\QQ{{\bf Q}}
\def\CC{{\bf C}}
\def\zmod#1{\,\,({\rm mod}\,\,#1)}

\font\twelvebf=cmbx12 at 15pt
{\eightpoint
\noindent Preliminary version \hfill Rome, January, 2008}
\hrule
\bigskip\bigskip
\centerline{\twelvebf A Modular Curve of Level 9 }
\bigskip
\centerline{\twelvebf and the Class Number One Problem}
\bigskip\bigskip
\centerline{\bf Burcu Baran}
\smallskip\centerline{\it Dipartimento di Matematica,}
\smallskip\centerline{\it Universit\`a di Roma ``Tor Vergata",}
\smallskip\centerline{\it I-00133  Roma, Italy.}
\bigskip\bigskip\bigskip

{\eightpoint
\noindent {\bf Abstract.} In this note we give an explicit parametrization of the modular curve associated to the normalizer of a non-split Cartan subgroup of level $9$. We determine all integral points of this modular curve. As an application, we give an alternative solution to the class number one problem.}

\beginsection 1.  Introduction

The class number one problem is the problem of determining exactly which imaginary quadratic fields have class number one. For a long time, it had been known that there are at least nine such imaginary quadratic fields. They are the ones with discriminant equal to
$-3$, $-4$, $-7$, $-8,$ $-11$, $-19$, $-43$, $-67$ and $-163$. It had long been conjectured that this is the complete list. The first proof was published in 1952 by Heegner $[{\bf{6}}]$. In his proof he used modular functions and reduced the issue to the problem of determining the integer solutions of certain diophantine equations. Unfortunately, he failed to prove some of the claims in his proof and so his proof was not accepted at the time. Later, in 1966, Baker proved in $[{\bf{1}}]$ that the list is complete. In his proof he used linear forms in logarithms. In the same year Stark also gave a proof in $[{\bf{13}}]$ using the modular functions  used by Heegner. Three years later, both Stark $[{\bf{14}}]$ and Deuring $[{\bf{4}}]$ independently patched up the supposed gap in Heegner's proof. In addition, in $1968$, Siegel $[{\bf{12}}]$ gave another solution to the class number one problem.
\smallskip
To every imaginary quadratic order $O$ of class number one, there is an associated elliptic curve $E_{O}$ over $\overline{\QQ}$ admitting complex multiplication by $O$. This curve is unique up to $\overline{\QQ}$-isomorphism and it can be defined over~$\QQ$. For any positive integer $n$, let $X_{ns}^{+}(n)$ denote the modular curve associated to the normalizer of a non-split Cartan subgroup of level $n$. It is defined over $\QQ$. If every prime $p$ that divides $n$ is inert in $O$, then $E_{O}$ gives rise to an integral point of~$X_{ns}^{+}(n)$; see $[{\bf{9}},\,{\rm{p}}.195]$. Here, by integral points, we mean points corresponding to elliptic curves with integral $j$-invariant. In this way, we can determine the imaginary quadratic fields with class number one by determining the integral points of $X_{ns}^{+}(n)$ for suitably chosen $n$. Some years ago, Serre pointed out $[{\bf{9}},\,{\rm{p}}.197]$ that, in fact, the solutions by Heegner and Stark of the class number one problem can be viewed in this way as the determination of the integral points of $X_{ns}^{+}(24)$. Later, following Serre's approach, Kenku gave a solution in $[{\bf{7}}]$ by using $X_{ns}^{+}(7)$ and Chen gave a moduli interpretation of Siegel's proof using the modular curve $X_{ns}^{+}(15)$ in $[{\bf{3}}]$. In this paper, we also follow Serre's approach and give a solution of the class number one problem by using the modular curve $X_{ns}^{+}(9)$. The following is the main result of this paper.
\smallskip\noindent
{\bf{Main result.}} We derive an explicit parametrization for the curve $X_{ns}^{+}(9)$ over $\QQ$. This modular curve is a genus zero curve with three cusps. Hence, by Siegel's theorem $[{\bf{9}},\,{\rm{p}}.95]$ it has finitely many integral points. We determine the integral points of this modular curve. In this way, we give another solution to the class number one problem.
\smallskip
Let $X(1)$ be the compactified modular curve whose affine part classifies isomorphism classes of elliptic curves. It is defined over $\QQ$. We choose the usual $j$-invariant as a uniformizer for $X(1)$. It is well-known  $[{\bf{3}}]$ that there exists a uniformizer $t$ defined over $\QQ$ for $X_{ns}^{+}(3)$ for which $j=t^{3}$.  In section $4$, we prove the following theorem.
\smallskip\noindent
\proclaim Theorem 1.1. There exist a uniformizer $y : X_{ns}^{+}(9) \longrightarrow {\bf{P}}^{1}$ defined over $\QQ$ such that the relation between $y$ and the uniformizer $t$ for $X_{ns}^{+}(3)$ that we mentioned above is given by
$$
t={{-3(y^{3}+3y^{2}-6y+4)(y^{3}+3y^{2}+3y+4)(5y^{3}-3y^{2}-3y+2)}\over{(y^{3}-3y+1)^{3}}}.
$$

\smallskip\noindent
The integral points of $X_{ns}^{+}(9)$ are the rational points for which $j$ and hence $t$ is in $\ZZ$. Thus, to find the integral points of $X_{ns}^{+}(9)$, we solve the equation in the above theorem for rational $y$ and integer $t$. In section $5$, we prove that such solutions correspond to the union of the integer solutions of the Diophantine equations $m^{3}-3mn^{2}+n^{3}=\pm1$ and $m^{3}-3mn^{2}+n^{3}=\pm3$. In the same section, we also give the integer solutions of these equations.
\smallskip
As a result, we find that there are nine integral points on the curve $X_{ns}^{+}(9)$. They are listed in Table $5.2$. One of these nine points correspond to an elliptic curve with $j$-invariant $1117947^{3}$. It  does not have complex multiplication. But, in some sense it behaves like it does modulo~$9$. The remaining eight points correspond to elliptic curves with complex multiplication by an imaginary quadratic order of class number one with the property that $3$ is inert in it. These are the ones with discriminant equal to $-4$, $-7$, $-16$, $-19$, $-28$, $-43$, $-67$, $-163$. Hence, this gives another proof that a tenth imaginary quadratic field with class number one does not exist. Indeed if such a field were to exist, then $3$ would be inert in that field, because in an imaginary quadratic field of class number one with discriminant $d$, all primes less than ${1+|d|}\over{4}$ are inert.  But then, it would give rise to an integral point on $X_{ns}^{+}(9)$, distinct from those that we have already found.

\beginsection 2. Non-split Cartan subgroups and their normalizers

In this section, we review non-split Cartan subgroups of GL$_{2}(\ZZ/n\ZZ)$ and their normalizers. Let $n$ be a positive integer and  $A$ be a finite free commutative ${\bf{Z}}/n{\bf{Z}}$-algebra of rank $2$ 
with unit discriminant. For example, if $R$ is an imaginary quadratic order with discriminant relatively prime to $n$, then $R/nR$ is such an algebra. By Galois theory, for a prime divisor $p$ of $n$, the $\FF_{p}$-algebra $A/pA$ is either ${\bf{F}}_{p}\times{\bf{F}}_{p}$, in which case $A$ is said to be {\it split} at $p$, or it is ${\bf{F}}_{p^2}$, in which case $A$ is said to be {\it non-split} at $p$. Fix $\{1, \alpha \}$ as a $\ZZ/n\ZZ$-basis for $A$. We have $A=(\ZZ/n\ZZ)[\alpha]$.  For any $a \in A^{\times}$, multiplication by $a$ is a $\ZZ/n\ZZ$-linear bijection. Therefore with respect to the $\ZZ/n\ZZ$-basis that we fixed, $A^{\times}$ embeds into GL$_{2}(\ZZ/n\ZZ)$.
\bigskip\noindent
{\bf{Definition 2.1.}} A  {\it{Cartan subgroup}} of ${\rm GL_{2}}(\ZZ/n\ZZ)$ is a subgroup
that arises as the image of such $A^{\times}$ in this way. If moreover $A$ is non-split at $p$ for every prime $p$ that divides~$n$, then the subgroup is 
called a {\it{non-split Cartan subgroup}} of ${\rm GL_{2}}(\ZZ/n\ZZ)$.
\bigskip\noindent
By Hensel's lemma, $A$ as above is unique up to isomorphism when it is non-split at all $p$ dividing $n$. Hence, all non-split Cartan subgroups of ${\rm GL_{2}}(\ZZ/n\ZZ)$ are conjugate to each other. Now, we describe the normalizer of a non-split Cartan subgroup of GL$_{2}(\ZZ/n\ZZ)$. Let $p$ be a prime that divides $n$ and $p^{r}$ be the maximum power of $p$ dividing $n$. There exists a unique ring automorphism $\sigma_{p}$ of order $2$ of $(\ZZ/n\ZZ)[\alpha]$ such that
 
 $$\eqalign{
 \sigma_{p}(\alpha) &\equiv \overline\alpha \ \ \ \ \ \ \ \ \ \ \,\,\,\,\,\,\,\,\,\,\pmod{p^{r}}, \cr
 \sigma_{p} &\equiv {\hbox{identity map}} \pmod{n/p^{r}}.\cr }
 $$
 where $\overline\alpha$ is the Galois conjugate of $\alpha$. Using the $\ZZ/n\ZZ$-basis $\{1, \alpha\}$ of A, we represent $\sigma_{p}$ by an element $S_{p}$ of ${\rm GL_{2}}(\ZZ/n\ZZ)$. The normalizer of a non-split Cartan subgroup $C$ of GL$_{2}(\ZZ/n\ZZ)$ is the group $$\langle\, C, S_{p} {\hbox{ for  }} p|n \,\rangle.$$
For the proof of this, see $[{\bf{2}}]$. The index of $C$ inside this subgroup is $2^{\nu}$ where $\nu$ is the number of prime divisors of $n$.
\bigskip\noindent
{\bf{Notation 2.2.}} We denote the intersection of the normalizer of a non-split Cartan subgroup of GL$_{2}(\ZZ/n\ZZ)$ with SL$_{2}(\ZZ/n\ZZ)$ by $C(n)$ and the subgroup of SL$_{2}(\ZZ)$ whose elements are congruent modulo $n$ to an element in $C(n)$ by $\Gamma_{C(n)}$.
\bigskip\noindent
In $[{\bf{2}}]$, we discuss non-split Cartan subgroups and their normalizers in general. In this paper, we mainly focus on the subgroups $C(3)$ and $C(9)$ and on the modular curves associated to these subgroups. Let $N$ be the kernel of the reduction modulo $3$ map 
$$
r:{\hbox{SL}}_{2}(\ZZ/9\ZZ) \longrightarrow {\hbox{SL}}_{2}(\ZZ/3\ZZ).
$$
This is an abelian group which is a vector space over $\FF_{3}$ of dimension $3$. The group $N$ is generated by the matrices 

$${\Big{\{}}\pmatrix{1&-3\cr 3&1\cr}, \pmatrix{-2&3\cr 3&4\cr}, \pmatrix{1&0\cr 3&1\cr}{\Big{\}}}.$$

The additive group of the ring $\ZZ/9\ZZ[i]:=\ZZ/9\ZZ[x]/(x^{2}+1)$ is a free rank $2$ module over $\ZZ/9\ZZ$. Choosing the basis $\{1, i\}$, we see that $C(9)$ is generated by $H$ and $\pmatrix{1&-3\cr 3&1\cr}$. Here, $H$ is the group generated by the matrices $\pmatrix{0&-1\cr 1&0\cr}$ and $\pmatrix{-1&-4\cr -4&1\cr}$. It is the $2$-Sylow subgroup of $C(9)$ and is isomorphic to the quaternion group of order $8$. By reducing $C(9)$ modulo $3$, we obtain the group $C(3)$. The group $H$ normalizes the subgroup $N'={\big{\langle}}\pmatrix{1&-3\cr 3&1\cr}{\big{\rangle}}$ of $N$ and it follows that we have the equality
$$
C(9)\,=\,N'H.$$
Of course, $H$ also normalizes $N$ itself and since $r$ induces an isomorphism from $H$ to $C(3)$, we also have the equality
$$
r^{-1}(C(3))\,=\,NH.$$
Direct computation shows that $H$ also normalizes the subgroup 
$$N''={\big{\langle}}\pmatrix{1&-3\cr 3&1\cr},\,\pmatrix{-2 &3\cr 3&4\cr}{\big{\rangle}},$$
where we have $N' \subset  N''  \subset N$. Therefore, we obtain the following inclusions of subgroups of $r^{-1}(C(3))$, each of index $3$,
$$
C(9)\,\, \subset \,\,N''H \,\, \subset \,\, r^{-1}(C(3)).$$
\bigskip\noindent
{\bf{Notation 2.3.}} We denote the subgroup $N''H$ of  $r^{-1}(C(3))$ by $B$.
\bigskip\noindent
As a consequence, we have the following commutative diagram,
$$
\def\normalbaselines{\baselineskip20pt\lineskip3pt
\lineskiplimit3pt}
\matrix{1& \longrightarrow &N& \longrightarrow &r^{-1}(C(3))&\mathop{\longrightarrow}\limits^{r}& C(3)& \longrightarrow &0 \cr
&&\bigcup &&\bigcup&&\parallel \cr
1& \longrightarrow &N''& \longrightarrow & B&\mathop{\longrightarrow}\limits^{r} & C(3)& \longrightarrow &0 \cr
&&\bigcup &&\bigcup&&\parallel \cr
1& \longrightarrow &N'& \longrightarrow &C(9)&\mathop{\longrightarrow}\limits^{r}& C(3)& \longrightarrow &0. \cr
}
$$\bigskip\noindent
{\bf{Notation 2.4.}} Let $\Gamma_{B}$ denote the subgroup of ${\hbox{SL}}_{2}(\ZZ)$ whose elements are congruent modulo $9$ to the elements in $B$. 
\bigskip\noindent
Then, we also have
$$
\Gamma_{C(9)} \,\, \subset \,\, \Gamma_{B} \,\, \subset  \,\, \Gamma_{C(3)}, 
$$
with $[\Gamma_{B}: \Gamma_{C(9)}]=3$ and $[\Gamma_{C(3)}:\Gamma_{B}]=3$.

\beginsection 3. The modular curves $X_{ns}^{+}(3)$, $X_{ns}^{+}(9)$ and $X_{B}$

In this section, we introduce the modular curves $X_{ns}^{+}(3)$, $X_{ns}^{+}(9)$ and $X_{B}$ over $\CC$. We also examine the natural covering maps between them.
\bigskip\noindent
{\bf{Definition 3.1.}} We define $$ \eqalign {X(1) &= {{\cal{H}^{*}}/ {{\hbox{SL}}_{2}(\ZZ)}} , \cr
X_{ns}^{+}(3) &= {{\cal{H}^{*}}/ {\Gamma_{C(3)}}} , \cr
X_{B} &= {{\cal{H}^{*}}/ {\Gamma_{B}}}, \cr
X_{ns}^{+}(9) &= {{\cal{H}^{*}}/ {\Gamma_{C(9)}}}. \cr }$$
Here, ${\cal{H}^{*}}={\cal{H}}\cup{\bf{P}}^{1}(\QQ)$ where $\cal{H}$ is the complex upper half plane. 
\bigskip\noindent
Consider the points $\infty$, $i$, $\rho=e^{{2\pi i}\over{3}}$ on ${\cal{H}^{*}}$. We have the natural covering maps
$$
X_{ns}^{+}(9) \,\, \mathop{\longrightarrow}\limits^{\pi_{1}} \,\, X_{B} \,\, \mathop{\longrightarrow}\limits^{\pi_{2}} \,\, X_{ns}^{+}(3) \,\, \mathop{\longrightarrow}\limits^{\pi_{3}} \,\, X(1).
$$
 The covering maps $\pi_{1}$, $\pi_{2}$ and $\pi_{3}$ all have degree $3$. The maps $\pi_{1}$, $\pi_{2}$ and $\pi_{3}$ ramify only above the points $\infty$, $\rho$, $i \in X(1)$. Thus, we study how the covering maps $\pi_{1}$, $\pi_{2}$ and $\pi_{3}$ branch above these points. This is the content of Proposition $3.3$ below. Before stating and proving this proposition, we recall the following facts which we will use in the proof.
 \smallskip
Let $\Gamma$ and $\Gamma'$ be subgroups of ${\hbox{SL}}_{2}(\ZZ)$ with $ \pm1 \in \Gamma \subset \Gamma'$. For every $\tau \in {\cal{H}^{*}}$, we define the stabilizer subgroups of $\tau$ by
$$
{\hbox{SL}}_{2}(\ZZ) _{\tau}= \{\gamma \in {\hbox{SL}}_{2}(\ZZ) : \, \, \gamma(\tau)=\tau \},
$$
$$
\Gamma _{\tau}= \{\gamma \in \Gamma : \, \, \gamma(\tau)=\tau \} \,\,\,\,\,\,\,\,\,\,{\hbox{and}}\,\,\,\,\,\,\,\,\,\,\Gamma' _{\tau}= \{\gamma \in \Gamma' : \, \, \gamma(\tau)=\tau \}.
$$
Consider the following commutative diagram,
$$
\def\normalbaselines{\baselineskip20pt\lineskip3pt
\lineskiplimit3pt}
\matrix{{\cal{H}^{*}}&=&{\cal{H}^{*}} \cr
\big\downarrow\rlap{$\scriptstyle \varphi'$} &&\big\downarrow\rlap{$\scriptstyle \varphi$} \cr
{\cal{H}^{*}} / {\Gamma}&\mathop{\longrightarrow}\limits^{f}&{\cal{H}^{*}} / \Gamma' \cr}
$$
where the vertical maps and $f$ are the natural projection maps. Let $z \in {\cal{H}^{*}}$ and $p=\varphi(z)$.  Let $q \in {\cal{H}^{*}} / {\Gamma}$ such that $f(q)=p$ and let $w \in {\cal{H}^{*}}$ be any point with $q=\varphi'(w)$. We have the following proposition.

\proclaim Proposition 3.2. The ramification index of $f$ at $q$ is equal to $[\Gamma'_{z} : \Gamma_{w}]$.

{\smallskip\noindent\bf{Proof:}} This is an easy special case of Proposition $1.37$ in $[{\bf{11}}]$.
\bigskip\noindent
In our case, the groups $\Gamma_{C(9)} \, \subset \, \Gamma_{B} \, \subset  \, \Gamma_{C(3)}$ contain $\pm1$. Now, we can prove the following proposition as promised above.

\proclaim Proposition 3.3. The covering maps $\pi_{1}$, $\pi_{2}$ and $ \pi_{3}$ branch above the points $\rho$, $i$, $\infty \in X(1)$ as in the following figure.
 \bigskip\medskip
 \centerline{\epsfxsize=5.5in \epsfbox[20 20 575 193]{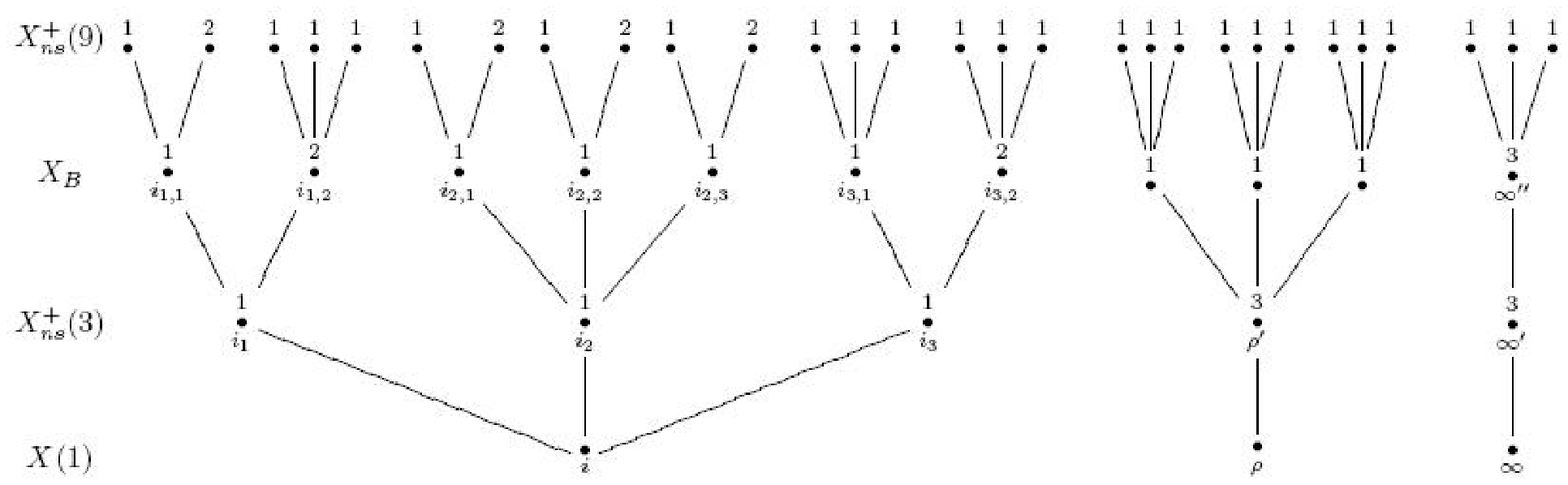}}
\bigskip\medskip\noindent
In the figure, the numbers over the points are the ramification indices of these points.

\smallskip\noindent
{\bf{Proof:}} First, we determine how the covering maps $\pi_{1}$, $\pi_{2}$ and $\pi_{3}$ branch above the point $\infty \in X(1)$. We know that SL$_{2}(\ZZ)_{\infty}$ is generated by the matrix $\pm\pmatrix{ 1 & 1 \cr 0&1\cr}$. Let $c$ be any cusp of $X$ for $X=X_{ns}^{+}(3)$ or $X_{B}$ or $X_{ns}^{+}(9)$. The matrices that are congruent to $\pm\pmatrix{ 1 & 0 \cr 0&1\cr}$ modulo $3$ are in $ \Gamma_{C(3)}$. Similarly, the matrices that are congruent to $\pm\pmatrix{ 1 & 0 \cr 0&1\cr}$ modulo $9$ are in $ \Gamma_{B}$ and also in $ \Gamma_{C(9)}$. Then we have,
$$
\Gamma_{c}={\pm\pmatrix{1&\ZZ \cr 0&1\cr}} \,\cap\, \Gamma= {\pm\pmatrix{1&a\ZZ \cr 0&1\cr}}
\,\,\,\,\,\,{\hbox{where}}\,\,\, {\cases { a = 3 & for $\Gamma=\Gamma_{C(3)}$, \cr
 a = 9 & for $\Gamma=\Gamma_{B}$ or $\Gamma_{C(9)}.$\cr}}$$
By Proposition $3.2$, the integer $a$ is the ramification index. Therefore, $X_{ns}^{+}(3)$ has one cusp with ramification index $3$ over $X(1)$, the curve $X_{B}$ has one cusp with ramification index $9$ over $X(1)$ and $X_{ns}^{+}(9)$ has three cusps with ramification indices $9$ over $X(1)$.
\smallskip
 Next, we determine how the covering maps $\pi_{1}$, $\pi_{2}$ and $ \pi_{3}$ branch above the point $\rho \in X(1)$. The group SL$_{2}(\ZZ)_{\rho}$ is generated by the matrix $\pmatrix{ 0 & -1 \cr 1&1\cr}$. Consider $\rho \in X_{ns}^{+}(3)$. The group $(\Gamma_{C(3)})_{\rho}$ is generated by the matrix $\pmatrix{ -1 & 0 \cr 0&-1\cr}$ which is the cube of the matrix $\pmatrix{ 0 & -1 \cr 1&1\cr}$. Thus, by Proposition $3.2$, above the point $\rho \in X(1)$, there exists only one point $\rho'$ with ramification index $3$ on the curve $X_{ns}^{+}(3)$. We know that a point that lies above the point $\rho \in X(1)$ can only have ramification index $1$ or $3$. Since $\rho' \in X_{ns}^{+}(3)$ has ramification index $3$, the curve $X_{B}$ has three points with ramification indices $3$ above the point $\rho \in X(1)$ and the curve $X_{ns}^{+}(9)$ has nine points with ramification indices $3$ above the point $\rho \in X(1)$.
 \smallskip
Finally, we determine how the covering maps $\pi_{1}$, $\pi_{2}$ and $\pi_{3}$ branch above the point $i \in X(1)$. First, we do this for the map $\pi_{3}$. The coset representatives of $\Gamma_{C(3)}$ in SL$_{2}(\ZZ)$ are $\pmatrix{ 1 & 0 \cr 0&1\cr}$, $\pmatrix{ 1 & 1 \cr 0&1\cr}$, $\pmatrix{ 1 & -1 \cr 0&1\cr}$. Thus, there exist at most three points $i$, $i+1$ and $i-1$ on $X_{ns}^{+}(3)$ over $i \in X(1)$. The groups SL$_{2}(\ZZ)_{i}$, SL$_{2}(\ZZ)_{i+1}$ and SL$_{2}(\ZZ)_{i-1}$ are generated by the matrices $\pmatrix{ 0 & -1 \cr 1&0\cr}$, $\pmatrix{ 1 & -2 \cr 1&-1\cr}$ and $\pmatrix{ -1 & -2 \cr 1&1\cr}$ respectively. These matrices are also elements of the group $\Gamma_{C(3)}$. Thus, by Proposition $3.2$, there are three points on $X_{ns}^{+}(3)$ with ramification indices $1$ over $i \in X(1)$. 
\smallskip
Next, we determine how the covering map $\pi_{2}$ branches above the points in $\pi_{3}^{-1}(i)$. The coset representatives of $\Gamma_{B}$ in $\Gamma_{C(3)}$ are $\pmatrix{ 1 & 0 \cr 0&1\cr}$, $\pmatrix{ 1 & 3 \cr 0&1\cr}$ and $\pmatrix{ 1 & -3 \cr 0&1\cr}$. This shows that there exist at most three points $i$, $i+3$ and $i-3$ on $X_{B}$ over $i \in X_{ns}^{+}(3)$. The matrix $\pmatrix{ 0 & -1 \cr 1&0\cr}$ which is the generator of the group $(\Gamma_{C(3)})_{i}$, is also an element of $\Gamma_{B}$. Thus, by Proposition $3.2$, the point $i \in X_{B}$ has ramification index $1$ over $i \in X_{ns}^{+}(3)$. A generator of the group $(\Gamma_{C(3)})_{i+3}$ is the product
$$
\pmatrix{ 1 & 3 \cr 0&1\cr}\,\pmatrix{ 0 & -1 \cr 1&0\cr}\,\pmatrix{ 1 & -3 \cr 0&1\cr},
$$
which is equal to $\pmatrix{ 3 & -1 \cr 1&-3\cr}$. This matrix is not inside the group $\Gamma_{B}$. We know that a point over $i \in X(1)$ can only have ramification index $1$ or $2$. Thus, the point  $i+3 \in X_{B}$ has ramification index $2$ over $i \in X_{ns}^{+}(3)$. But, the degree of $\pi_{2}$ is $3$. This shows that the points $i+3$ and $i-3$ are equal on $X_{B}$. As a result, there exist two points $i$, $i+3$ over $i \in X_{ns}^{+}(3)$ with ramification indices $1$, $2$ respectively. Similar calculations show that there exist three points $i+1$, $i+4$ and $i-2$ on $X_{B}$ over the point $i+1 \in X_{ns}^{+}(3)$ all with ramification indices $1$  and there exist two points $i-1$, $i-4$ on $X_{B}$ over the point $i-1 \in X_{ns}^{+}(3)$ with ramification indices $2$, $1$ respectively.
\smallskip
 Finally, we determine how the covering map $\pi_{1}$ branches at the points lying over $i$. Since, the points $i+3$ and $i-1 \in X_{B}$ are already ramified over the point $i \in X(1)$, they are unramified for the map $\pi_{1}$. Thus, we only need to consider the points $i$, $i+1$, $i+4$, $i-2$ and $i-4 \in X_{B}$. Coset representatives of $\Gamma_{C(9)}$ inside $\Gamma_{B}$ are $\pmatrix{ -2 & 3 \cr 3& 4\cr}$, $\pmatrix{ 4 & -3 \cr -3&-2\cr}$ and $\pmatrix{ 1 & 0 \cr 0&1\cr}$, so there exist at most three points on $X_{ns}^{+}(9)$ over $i \in X_{B}$. One of these three points is $i$. The generator of $(\Gamma_{B})_{i}$ is the matrix $\pmatrix{ 0 & -1 \cr 1&0\cr}$ which is an element of $\Gamma_{C(9)}$. Thus, by Proposition $3.2$, the point $i \in X_{ns}^{+}(9)$ has ramification index $1$ over $i \in X_{B}$. The stabilizer in $\Gamma_{B}$ of one of the other two points over $i \in X_{B}$ is generated by the product of the matrices
$$\pmatrix{ -2 & 3 \cr 3& 4\cr}\,\pmatrix{ 0 & -1 \cr 1& 0\cr}\,\pmatrix{ 4 & -3 \cr -3& -2\cr}$$
which is equal to $\pmatrix{ -3 & -4 \cr -2& 3\cr}$. This matrix is not an element of $\Gamma_{C(9)}$. Thus, this point has ramification index $2$ over $i \in X_{B}$. This shows that these two points are equal on~$X_{ns}^{+}(9)$. Consequently, there exist two points over $i \in X_{B}$ with ramification indices $1$ and~$2$. Similar calculations show that there exist two points on $X_{ns}^{+}(9)$ over $i+1$, $i+4$ and $i-2 \in X_{B}$ with ramification indices 1 and 2 and there exist three points on $X_{ns}^{+}(9)$ over $i-4 \in X_{B}$ all with ramification indices $1$. Hence, we proved the proposition.

\proclaim Corollary 3.4. The modular curves $X_{ns}^{+}(3)$, $X_{B}$ and $X_{ns}^{+}(9)$ all have genus $0$.
 
 {\smallskip\noindent\bf{Proof:}} Proposition $3.3$ says how the covering maps $\pi_{1}$, $\pi_{2}$ and $\pi_{3}$ branch above the points $i$, $\rho$, $\infty \in X(1)$. These covering maps are unramified above every other point of $X(1)$. Thus, the corollary follows from the Riemann-Hurwitz formula $[{\bf{5}},\,{\rm{p}}. 66]$.
\bigskip\noindent

For a formula for the genus of the modular curve $X_{ns}^{+}(n)$ where $n$ is any positive integer, see $[{\bf{2}}]$.

\beginsection 4. An explicit parametrization over $\QQ$ for $X_{ns}^{+}(9)$ over $X(1)$ 

 \noindent
In this section, following Chen's method $[{\bf{3}}]$, we obtain an explicit parametrization over $\QQ$ for the modular curve $X_{ns}^{+}(9)$ over~$X(1)$. In the previous section, we obtained enough information about the covering maps $\pi_{1}$, $\pi_{2}$ and $\pi_{3}$ to write down such a parametrization. 

\bigskip\noindent
{\bf Definition 4.1.}  Let $X$ be a projective non-singular curve defined over $\QQ$. Suppose it has genus zero and at least one $\QQ$-rational point. Then, there exists an isomorphism 
$$
h: X \longrightarrow {\bf{P}}^{1},
$$
which is  defined over $\QQ$. It is unique up to an automorphism of ${\bf{P}}^{1}_{\QQ}$. We call the map $h$ a {\it uniformizer} over $\QQ$ for $X$.
\bigskip\noindent
Consider the modular curve $X(1)$ and suppose there is given a finite covering map
$$
\upsilon : X \longrightarrow X(1),
$$
over $\QQ$. We have the induced map $\upsilon^{*}$ on the function fields of these curves. Let $j$ be the unique uniformizer over $\QQ$ of $X(1)$ which is the usual $j$-function. Identify $\upsilon^{*}(j)$ with~$j$. Then, by an {\it{explicit parametrization}} over $\QQ$ for $X$ over $X(1)$ with respect to the uniformizers $h$ and $j$, we mean the unique relation
$$
j=\lambda {{P(h)}\over{Q(h)}},
$$
where $\lambda \in \QQ^{\times}$ and $P$, $Q$ are monic coprime polynomials with rational coefficients. These notions also make sense over any extension of $\QQ$. Working over $\CC$, we see that the constant $\lambda$ and the polynomials $P$, $Q$ are given by
$$\eqalign{
\lambda &= 1728 {{Q(h(z_{0}))}\over{P(h(z_{0}))}}, \cr
P(T) &= \prod_{z \in \upsilon^{-1}(\rho)}(T- h(z))^{e(z)}, \cr
Q(T) &= \prod_{z \in \upsilon^{-1}(\infty)}(T- h(z))^{e(z)}, \cr}
$$
where $e(z)$ is the ramification index of $\upsilon$ at $z \in X$ and $z_{0}$ is any point in $\upsilon^{-1}(i)$.
\bigskip\noindent
{\bf Remark 4.2.} Since the determinant map from the normalizer of a non-split Cartan subgroup of ${\hbox{GL}}_{2}(\ZZ/n\ZZ)$ to $(\ZZ/n\ZZ)^{\times}$ is surjective, by $[{\bf{ 5}}, {\hbox {Thm}}\,7.6.3]$, $[{\bf{9}}, {\hbox{p}}.\,194]$ the modular curve $X_{ns}^{+}(n)$ can be defined over $\QQ$. In this paper, we use this natural choice of $\QQ$-structure on the complex modular curves $X_{ns}^{+}(3)$ and $X_{ns}^{+}(9)$, so $\pi_{3}$ and $\pi_{2} \circ \pi_{1}$ are defined over $\QQ$.

\medskip\noindent
{\bf Remark 4.3.} The group $B$ is the intersection of ${\hbox {SL}}_2(\ZZ/9\ZZ)$ and the subgroup $B'$  of ${\hbox {GL}}_2(\ZZ/9\ZZ)$ generated by $B$ and the scalar matrices. Since the image of $B'$ under the determinant map is the subgroup of squares of $(\ZZ/9\ZZ)^{\times}$, a refinement of the determinant arguments in $[{\bf {5}}, {\hbox {section}}\,7.6]$, $[{\bf{9}}, {\hbox{p}}.\,194]$ shows that the modular curve $X_{B}$ is (uniquely) defined over $\QQ(\sqrt{-3})$ compatibly with the $\QQ$-structure on $X_{ns}^{+}(3)$ and $X_{ns}^{+}(9)$. By Lemma~$4.5.$ below, it can not be defined over $\QQ$. In this paper, we use this natural choice of $\QQ(\sqrt{-3})$-structure on the complex modular curve $X_{B}$.
\bigskip
Now first, we will write an explicit parametrization over $\QQ$ for $X_{ns}^{+}(3)$ over $X(1)$, then over $\QQ(\sqrt{-3})$ for $X_{B}$ over $X_{ns}^{+}(3)$ and then again over $\QQ(\sqrt{-3})$ for $X_{ns}^{+}(9)$ over~$X_{B}$. In the end, we will obtain an explicit parametrization over $\QQ(\sqrt{-3})$ for the curve $X_{ns}^{+}(9)$ over $X(1)$ and finally we will compute Gal$(\QQ(\sqrt{-3})/\QQ)$-actions to find an explicit parametrization over $\QQ$ for the curve $X_{ns}^{+}(9)$ over $X(1)$.
\smallskip
We start by obtaining an explicit parametrization over $\QQ$ of the modular curve $X_{ns}^{+}(3)$ over~$X(1)$. In fact, this is well-known $[{\bf{3}}]$, but we recall the argument for the convenience of the reader. Consider the natural covering map 
$$
\pi_{3} : X_{ns}^{+}(3) \longrightarrow X(1),
$$
over $\QQ$. As we have shown in Proposition $3.3$, there exists a unique point $\rho' \in X_{ns}^{+}(3)$ lying over $\rho \in X(1)$ and a unique point $\infty' \in X_{ns}^{+}(3)$ lying above $\infty \in X(1)$. Hence, by uniqueness, these are $\QQ$-points of $X_{ns}^{+}(3)$. By Corollary~$3.4$, the curve $X_{ns}^{+}(3)$ has genus zero. Thus, we choose the unique uniformizer
$$
t : X_{ns}^{+}(3) \longrightarrow {\bf{P}}^{1},
$$
over $\QQ$, such that $t(\rho')=0$ and $t(\infty')=\infty$. As the degree-$3$ map $\pi_{3}$ is totally ramified over $0$ and $\infty$ and nowhere else, the relation between $j$ and $t$ has the form   $j=\lambda t^{3}$ where $\lambda \in \QQ^{\times}$ is determined up to a cube in $\QQ^{\times}$. Let $E_{\ZZ[i]}$ denote the unique elliptic curve over $\CC$, up to isomorphism, with endomorphism ring $\ZZ[i]$. The $j$-invariant of $E_{\ZZ[i]}$ is $1728=12^{3}$. Since $3$ is inert in $\ZZ[i]$, the moduli interpretation of $X_{ns}^{+}(3)$ shows that $t(E_{\ZZ[i]}) \in \QQ$. Thus, the constant $\lambda$ is a cube in $\QQ^{\times}$. Hence, we can assume that $\lambda$ is equal to $1$ and we have 
$$
j=t^{3}. \eqno{(4.1)}
$$
Next, we obtain an explicit parametrization over $\QQ(\sqrt{-3})$ for the curve $X_{B}$ over $X_{ns}^{+}(3)$. For the curve $X_{ns}^{+}(3)$, we use the uniformizer $t$ over $\QQ$ that we fixed above. Consider the natural covering map
$$
\pi_{2} : X_{B} \longrightarrow X_{ns}^{+}(3),
$$
over $\QQ(\sqrt{-3})$. In Proposition $3.3$, we have determined how this  covering map viewed over $\CC$ is branched above $i$, $\rho$, $\infty \in X(1)$. Note that from equation $(4.1)$, for the points $i_{1}$, $i_{2}$, $i_{3} \in X_{ns}^{+}(3)$ over $i \in X(1)$, the values of $t(i_{1})$, $t(i_{2})$, $t(i_{3})$ are the roots of the polynomial 
$$T^{3}-j(i)=T^{3}-1728=T^{3}-12^{3}.$$
Since $t(E_{\ZZ[i]}) \in \QQ$, the point $i_{1}=i \in X_{ns}^{+}(3)$ is $\QQ$-rational. Thus, we have $t(i_{1})=12$ and we can uniquely choose the labels $i_{2}$ and $i_{3}$ so that $t(i_{3})=12\rho$ and $t(i_{2})=12\rho^{-1}$. As we mentioned before, $\rho={{-1+{\sqrt{-3}}} \over 2}$ is a primitive third root of unity.

\proclaim Proposition 4.4. Using the notation of Proposition $3.3$, there exists a unique uniformizer 
$$
w : X_{B} \longrightarrow {\bf{P}}^{1}
$$
over $\QQ(\sqrt{-3})$, such that $w(\infty'')=\infty$, $w(i_{3,1})=2{\sqrt{-3}}$ and $w(i_{3,2})=-{\sqrt{-3}}$. The points $i_{3,1}$ and $i_{3,2}$ are the unique points lying over $i_{3} \in X_{ns}^{+}(3)$, for which $t(i_{3})=12\rho$, with ramification indices $1$ and $2$ respectively. Here, $t$ is the uniformizer of $X_{ns}^{+}(3)$ that we fixed above. The uniformizers $t$ and $w$ are related by 
$$
t=\rho^{-1}(w^{3}+9w-6),
$$
where $\rho$ is the third root of unity ${-1+{\sqrt{-3}}}\over 2$.

\smallskip\noindent
{\bf{Proof:}} Consider the figure in Proposition $3.3$. We see that $\infty'' \in X_{B}$ is the unique point above $\infty' \in X_{ns}^{+}(3)$, the point $i_{3,1} \in X_{B}$ is the unique unramified point above $i_{3} \in X_{ns}^{+}(3)$ and $i_{3,2} \in X_{B}$ is the unique ramified point above $i_{3} \in X_{ns}^{+}(3)$. By uniqueness, these points are $\QQ(\sqrt{-3})$-rational. Therefore, we can choose the uniformizer
$$
\eta : X_{B} \longrightarrow {\bf{P}}^{1}
$$
over $\QQ(\sqrt{-3})$ such that $\eta(\infty'')=\infty$, $\eta(i_{3,1})=1$ and $\eta(i_{3,2})=0$. Then the relation between $t$ and $\eta$ has the form
$$
t=\lambda{\prod_{i=1}^{3}}(\eta-\eta(\rho_{i}))=\lambda(\eta^{3}+A\eta^{2}+B\eta+C),
$$
where $\lambda \in \QQ(\sqrt{-3})^{\times}$ and $A$, $B$, $C \in \QQ(\sqrt{-3})$. Evaluating this equality at $i_{3, 2}$ yields
$$
\lambda C = t(i_{3})=12 \rho.
$$
But $\pi_{2}$ has ramification index $2$ at $i_{3, 2}$, so $t-t(i_{3})$ has double zero at $i_{3, 2}$. Since we have
$$
 t-t(i_{3})=t-\lambda C=\lambda(\eta^{3}+A\eta^{2}+B\eta),
 $$
 we conclude that $\lambda B=0$, so $B=0$.
Thus, we get
$$
t=\lambda(\eta^{3}+A\eta^{2}+C)
$$
where $\lambda C=12\rho$. We also have $\eta(i_{3, 1})=1$, so we get 
$$
\lambda(1+A+C)=t(i_{3})=\lambda C.
$$
This shows that $A=-1$ and $t=\lambda(\eta^{3}-\eta^{2}+C)$.
\smallskip
If we similarly consider the way that the covering map $\pi_{2}$ branches above the point $i_{1} \in X_{ns}^{+}(3)$, we see that $t=\lambda(\eta-D)(\eta-E)^{2}+12$ in $\QQ(\sqrt{-3})[\eta]$ where $D=\eta(i_{1,1})$ and $E=\eta(i_{1,2})$. Computing coefficients in our two expressions for $t$ yields the system of equations
$$\eqalign  {
\lambda(-1+2E+D)&=0, \cr
\lambda(E^{2}+2DE)&=0, \cr
\lambda C-12+\lambda DE^{2}&=0. \cr}
$$
Since $\lambda$ is non-zero, the first equation says $D=1-2E$, so the second equation yields $E=0$ or $E=2/3$. But $E=0$ gives a contradiction with the third equation, since we know that $\lambda C=12\rho$. Therefore, $E=2/3$ and hence $D=-1/3$. Using the equality $\lambda C=12\rho$, the third equation gives that $\lambda=81(\rho-1)$ and $C={{-4(\rho-1)}\over 81}$. Hence, the relation between $t$ and $\eta$ is 
$$
t=-81(\rho-1)(\eta^{3}-\eta^{2}+{{-4(\rho-1)}\over 81}).
$$
We can simplify the equation for $t$ in terms of $\eta$  by a change of variable. We put $w=-(\sqrt{-3})^{3}\eta-{\sqrt{-3}}$ to get the equation
$$
t=\rho^{-1}(w^{3}+9w-6).
$$
With this change of variable, we also have $w(\infty'')=\infty$, $w(i_{3,1})=2{\sqrt{-3}}$ and~$w(i_{3,2})=-{\sqrt{-3}}$. Hence, the proposition follows.  
\bigskip\noindent
Next, our aim is to find the relation over $\QQ$ between $w$ and a uniformizer of the modular curve~$X_{ns}^{+}(9)$. Extending the constant field of the function fields  $\QQ(j)$, $\QQ(t)$, $\QQ(X_{ns}^{+}(9))$  of the curves $X(1)$, $X_{ns}^{+}(3)$ and $X_{ns}^{+}(9)$ respectively by $\rho$, we have the following diagram relating the curves and their function fields over $\QQ(\sqrt{-3})$.
$$
\def\normalbaselines{\baselineskip20pt\lineskip3pt
\lineskiplimit3pt}
\matrix{X_{ns}^{+}(9)& \longleftrightarrow &\QQ(\sqrt{-3})(X_{ns}^{+}(9)) \cr
\big\downarrow\rlap{$\scriptstyle \pi_{1}$} && \cup \cr
X_{B}& \longleftrightarrow &\QQ(\sqrt{-3}, w) \cr
\big\downarrow\rlap{$\scriptstyle \pi_{2}$} && \cup \cr
X_{ns}^{+}(3)& \longleftrightarrow &\QQ(\sqrt{-3}, t) \cr
\big\downarrow\rlap{$\scriptstyle \pi_{3}$} &&\cup \cr
X(1)& \longleftrightarrow &\QQ(\sqrt{-3}, j) \cr
}
$$
The key point is the following lemma.

\proclaim Lemma 4.5. The intermediate curve $X_{B}$ which is defined over $\QQ(\sqrt{-3})$ (by Remark $4.3$) is {\it{not}} defined over $\QQ$ with respect to the natural $\QQ$-structure on $X_{ns}^{+}(9)$.

\smallskip\noindent
{\bf{Proof:}} Suppose $X_{B}$ is defined over $\QQ$ compatibly with $X_{ns}^{+}(9)$. Then, the covering map 
$$
\pi_{2} : X_{B} \longrightarrow X_{ns}^{+}(3).
$$
is defined over $\QQ$, because $\pi_{2} \circ \pi_{1}$ is defined over $\QQ$. Consider the figure in Proposition~$3.3$. The map $\pi_{2}$ is ramified over $i_{3}$, but not over its $\QQ$-conjugate point $i_{2}$. Thus, $\pi_{2}$ is not defined over $\QQ$. Therefore, the curve $X_{B}$ is not defined over $\QQ$, as required.
\bigskip\noindent
Let $\sigma$ be a generator of the group ${\hbox{Gal}}(\QQ(\sqrt{-3})/\QQ)$. The action of $\sigma$ naturally extends to an action on $\QQ(\sqrt{-3})(X_{ns}^{+}(9))$. From Lemma $4.5$, it follows that $\sigma$ cannot preserve the subfield $\QQ(\sqrt{-3})(X_{B})=\QQ(\sqrt{-3}, w)$. Thus, $\sigma$ carries $\QQ(\sqrt{-3}, w)$ to some other subfield $\QQ(\sqrt{-3}, w')$ over $\QQ(\sqrt{-3})(X_{ns}^{+}(3))$ where $w'=\sigma(w)$.
\smallskip
To summarize, we have the following diagram of degree $3$ extensions.
$$
\def\normalbaselines{\baselineskip20pt\lineskip3pt
\lineskiplimit3pt}
\matrix{&&\QQ(\sqrt{-3})(X_{ns}^{+}(9))&& \cr
&\nearrow&&\nwarrow& \cr
\QQ(\sqrt{-3}, w)&&&&\QQ(\sqrt{-3}, w') \cr
&\nwarrow&&\nearrow& \cr
&&\QQ(\sqrt{-3},t)&& \cr
}
$$
Here, the field extension $\QQ(\sqrt{-3}, w')$ over $\QQ(\sqrt{-3}, t)$ is given by the equation
$$
t= \rho(w'^{3}+9w'-6).
$$
The field $\QQ(\sqrt{-3}, w, w')$ lies between the fields $\QQ(\sqrt{-3})(X_{ns}^{+}(9))$ and $\QQ(\sqrt{-3}, w)$. Since the degree of the extension $\QQ(\sqrt{-3})(X_{ns}^{+}(9))/\QQ(\sqrt{-3}, t)$ is $9$ and the field $\QQ(\sqrt{-3}, w)$ is not equal to the field $\QQ(\sqrt{-3}, w')$, we have
$$
\QQ(\sqrt{-3})(X_{ns}^{+}(9))=\QQ(\sqrt{-3}, w, w').$$
The elements $w$ and $w'$ satisfy
$$ 
 \rho^{-1}(w^{3}+9w-6)= \rho(w'^{3}+9w'-6). \eqno{(4.2)}
$$
It can be shown that equation $(4.2)$ defines an irreducible cubic curve with unique singular point $(w, w')=({\sqrt{-3}}, -{\sqrt{-3}})$, and hence is a singular model for $X_{ns}^{+}(9)$ over $\QQ(\sqrt{-3})$. But we do not need to know this. However, it motivates the definition of the ``slope parameter" $u$ in the proof of Proposition $4.6$. In the following proposition, we find a uniformizer defined over $\QQ$ for the curve $X_{ns}^{+}(9)$.

\proclaim Proposition 4.6. There exists a uniformizer $y : X_{ns}^{+}(9) \longrightarrow {\bf{P}}^{1}$ over $\QQ$ such that the relation between $y$ defined over $\QQ$ and the uniformizer $t$ defined over $\QQ$ of $X_{ns}^{+}(3)$ that we fixed above is 
$$
t={{-3(y^{3}+3y^{2}-6y+4)(y^{3}+3y^{2}+3y+4)(5y^{3}-3y^{2}-3y+2)}\over{(y^{3}-3y+1)^{3}}}.
$$

\smallskip\noindent
{\bf{Proof:}} We follow the above notation. Consider the function
$$
u={{w-{\sqrt{-3}}} \over {w'+{\sqrt{-3}}}},
$$
in $\QQ(\sqrt{-3},w,w')$. Recall that we have $\sigma(w)=w'$ where $\sigma$ is the generator of the group Gal$(\QQ(\sqrt{-3})/\QQ)$. It follows that $\sigma(u)={1\over u}$. We change variables 
$$w-{\sqrt{-3}}=\vartheta\,\,\,\,\,\,\,\,\,\,\,\,\,\,\,{\hbox{and}}\,\,\,\,\,\,\,\,\,\,\,\,\,\,\,w'+{\sqrt{-3}}=\tau.  \eqno{(4.3)}$$
Inserting them into equation $(4.2)$, we  obtain
$$
(\vartheta+{\sqrt{-3}})^{3}+9(\vartheta+{\sqrt{-3}})-6=\rho^{2}((\tau-{\sqrt{-3}})^{3}+9(\tau-{\sqrt{-3}})-6). 
$$
We simplify this equation and we get
$$
\vartheta^{3}+3{\sqrt{-3}}\vartheta^{2}=\rho^{-1}\tau^{3}+\big({{-9+3{\sqrt{-3}}}\over2}\big)\tau^{2}. \eqno{(4.4)}
$$

Substituting $\vartheta=u \tau$ in $(4.4)$ and inserting the solutions into equations $(4.3)$, we can write $w'$ and $w$ in terms of $u$. For $w$, we get
$$
w=3u{\sqrt{-3}}\big({{-u^{2}-\rho^{-1}}\over {u^{3}-\rho^{-1}}}\big)+{\sqrt{-3}}.
$$
Since $u$ is a uniformizer defined over $\QQ(\sqrt{-3})$ for the curve $X_{ns}^{+}(9)$, this is the relation between the uniformizer $u$ of the curve $X_{ns}^{+}(9)$ and the uniformizer $w$ of the curve $X_{B}$. As we also have $t=\rho^{-1}(w^{3}+9w-6)$ where $t$ is the uniformizer of the curve $X_{ns}^{+}(3)$, we obtain
$$
t=\rho^{-1}\big(\big(3{\sqrt{-3}}u\big({{-u^{2}-\rho^{-1}}\over {u^{3}-\rho^{-1}}}\big)+{\sqrt{-3}}\big)^{3}+9\big(3{\sqrt{-3}}u\big({{-u^{2}-\rho^{-1}}\over {u^{3}-\rho^{-1}}}\big)+{\sqrt{-3}}\big)-6\big). \eqno{(4.5)}
$$
But this equation has coefficients in $\QQ(\sqrt{-3})$. As we know, the curves $X_{ns}^{+}(9)$, $X_{ns}^{+}(3)$ and their projections to $X(1)$ are defined over $\QQ$. Thus, there exists a uniformizer $y$ defined over $\QQ$ for $X_{ns}^{+}(9)$ such that $t \in \QQ(y)$. We now find such $y$.  As we showed above, $\sigma(u)={1\over u}$. Thus, if we make the change of variable
$$
u={{y+\rho}\over{\rho y+1}}
$$
then we find $\sigma(y)=y$. Hence, $\QQ(X_{ns}^{+}(9))=\QQ(y)$. We also compute
$$
t={{-3(y^{3}+3y^{2}-6y+4)(y^{3}+3y^{2}+3y+4)(5y^{3}-3y^{2}-3y+2)}\over{(y^{3}-3y+1)^{3}}},
$$
using $(4.5)$. Therefore, the proposition follows.
\bigskip\noindent
Combining equation $(4.1)$ and the equation in Proposition $4.6$ gives an explicit parametrization over $\QQ$ for $X_{ns}^{+}(9)$ over $X(1)$ in the sense defined early in this section.

\beginsection 5. The integral points and imaginary quadratic fields of class number one

\noindent
In this section, we determine the integral points on the curve $X_{ns}^{+}(9)$ and obtain the imaginary quadratic fields of class number one. As we mentioned before, by integral points on the curve $X_{ns}^{+}(9)$, we mean the points corresponding to curves with integral $j$-invariant, i.e., the rational points on $X_{ns}^{+}(9)$ whose image under the covering map $\pi_{1} \circ \pi_{2} \circ \pi_{3}$ of section~$3$ are integral points on $X(1)$. 
\smallskip
Consider equation $(4.1)$ and the equation in Proposition~$4.6$. We see that, to find integral points on $X_{ns}^{+}(9)$, we should search for the points $(y,\,t)$ which satisfy the equation in Proposition $4.6$ for $t \in \ZZ$ and $y \in \QQ \cup \{ \infty \}$. Let $y={m\over n}$ where $m, n \in \ZZ$ with ${\hbox{gcd}}(m, n)=1$. We have
$$
\eqalignno{ &t\,\,=\,\,&(5.1)\cr
&{{-3(m^{3}+3m^{2}n-6mn^{2}+4n^{3})(m^{3}+3m^{2}n+3mn^{2}+4n^{3})(5m^{3}-3m^{2}n-3mn^{2}+2n^{3})}\over{(m^{3}-3mn^{2}+n^{3})^{3}}}. \cr}
$$

\proclaim Lemma 5.1. For any integer solution $(m, n, t)$ of equation $(5.1)$ with gcd$(m, n)=1$, we have either $m^{3}-3mn^{2}+n^{3}=\pm1$ or $m^{3}-3mn^{2}+n^{3}=\pm3$.

\smallskip\noindent
{\bf{Proof:}} Suppose $(m, n, t)$ is an integer solution to $(5.1)$ with gcd$(m, n)=1$. If $n$ is equal to $0$ then the only solution is $(1,\,0,\,-15)$. Thus, to see the other solutions, we assume that $n$ is not equal to $0$. Suppose $$m^{3}-3mn^{2}+n^{3} \equiv 0 \zmod p$$ 
for some prime $p$. Note that if $n \equiv 0 \zmod p$ then $m \equiv 0 \zmod p$. But this contradicts the assumption gcd$(m, n)=1$, so $n \not\equiv 0 \zmod p$. Thus, we can say that 
$$
c^{3}-3c+1 \equiv 0 \zmod p
$$
where $c={m \over n}$. Since $t$ is an integer we therefore have,
$$-3(c^{3}+3c^{2}-6c+4)\,(c^{3}+3c^{2}+3c+4)\,(5c^{3}-3c^{2}-3c+2) \equiv 0 \zmod p. $$ As we also have 
$${\hbox{resultant}}((y^{3}+3y^{2}-6y+4)(y^{3}+3y^{2}+3y+4)(5y^{3}-3y^{2}-3y+2),\,y^{3}-3y+1)=3^{15},$$ it follows that $p=3$. But, the equation $y^{3}-3y+1 \equiv 0 \zmod {3^{2}}$ does not admit any solution $y \in \ZZ/9\ZZ$. Therefore, we have $m^{3}-3mn^{2}+n^{3}=\pm3$ or $\pm 1$ and the lemma follows.
\bigskip\noindent
Hence, to find all integer solutions $(m, n, t)$ of equation $(5.1)$ with gcd$(m, n)=1$, it is enough to solve the two cubic Diophantine equations,
$$
m^{3}-3mn^{2}+n^{3}=\pm1\,\,\,\,\,\,\,\,\,\,\,\,\,{\hbox{and}}\,\,\,\,\,\,\,\,\,\,\,\,\, m^{3}-3mn^{2}+n^{3}=\pm3.
$$
Ljunggren $[{\bf{8}}]$, shows that the only integer solutions of the equation $m^{3}-3mn^{2}+n^{3}=1$ are $(2,-1)$, $(-3,-2)$, $(-1,-1)$, $(1,0)$, $(1,3)$, $(0,1)$. He uses Skolem's method to determine these solutions. It can be shown in the same way, again by using Skolem's method, that the only integer solutions of the equation $m^{3}-3mn^{2}+n^{3}=3$ are $(-1,-2)$, $(-1,1)$ and $(2,1)$. Alternatively, these equations are Thue equations and the PARI/GP computer package solves these Thue equations in less than one second.
\smallskip
As a result, we see that there are nine integral points on the modular curve $X_{ns}^{+}(9)$. By using equations $(4.1)$ and $(5.1)$ we calculate the $j$-invariants of the corresponding elliptic curves. Also, by using the table in $[{\bf{9}},\,{\rm{p}}.\,192]$ we determine the discriminants $d$ of the corresponding imaginary quadratic orders $R_{d}$ of class number one. Table $5.2$ shows all of these data.

\bigskip\medskip
\vbox{\noindent {\bf Table 5.2.}
\bigskip
\centerline {\vbox {\offinterlineskip
\hrule\halign{&\vrule#&\strut\quad\hfil#\quad\cr
height3pt
&\omit&&\omit&&\omit&\cr
&integral solution $(m,n)$&&$j$-invariant &&discriminant $d$&\cr
height3pt
&\omit&&\omit&&\omit&\cr
\noalign{\hrule}
height2pt
&\omit&&\omit&&\omit&\cr
&$(-1,1)$&&$2^{{\hbox{$ \scriptstyle6$}}}3^{{\hbox{$ \scriptstyle3$}}}$&&$-4$&\cr
&(1,0)&&$-3^{{\hbox{$ \scriptstyle3$}}}5^{{\hbox{$ \scriptstyle3$}}}$&
&$-7$&\cr
&$(-1,-1)$&&$2^{{\hbox{$ \scriptstyle3$}}}3^{{\hbox{$ \scriptstyle3$}}}11^{{\hbox{$ \scriptstyle3$}}}$&&$-16$&\cr
&(0,1)&&$-2^{{\hbox{$ \scriptstyle15$}}}3^{{\hbox{$ \scriptstyle3$}}}$&&$-19$&\cr
&$(-1,-2)$&&$3^{{\hbox{$ \scriptstyle3$}}}5^{{\hbox{$ \scriptstyle3$}}}17^{{\hbox{$ \scriptstyle3$}}}$&
&$-28$&\cr
&(2,1)&&$-2^{{\hbox{$ \scriptstyle18$}}}3^{{\hbox{$ \scriptstyle3$}}}5^{{\hbox{$ \scriptstyle3$}}}$&&$-43$& \cr
&$(2,-1)$&&$-2^{{\hbox{$ \scriptstyle15$}}}3^{{\hbox{$ \scriptstyle3$}}}5^{{\hbox{$ \scriptstyle3$}}}11^{{\hbox{$ \scriptstyle3$}}}$&&$-67$&\cr
&(1,3)&&$-2^{{\hbox{$ \scriptstyle18$}}}3^{{\hbox{$ \scriptstyle3$}}}5^{{\hbox{$ \scriptstyle3$}}}23^{{\hbox{$ \scriptstyle3$}}}29^{{\hbox{$ \scriptstyle3$}}}$&&$-163$&\cr
&$(-3,-2)$&&$3^{{\hbox{$ \scriptstyle3$}}}41^{{\hbox{$ \scriptstyle3$}}}61^{{\hbox{$ \scriptstyle3$}}}149^{{\hbox{$ \scriptstyle3$}}}$&&&\cr
height3pt
&\omit&&\omit&&\omit&\cr
}\hrule}}}
\bigskip\medskip
\noindent
In Table $5.2$, we find the imaginary quadratic fields of class number one with discriminant equal to $-4$, $-7$, $-19$, $-43$, $-67$ and $-163$. The imaginary quadratic fields  of class number one with discriminant equal to $-3$, $-8$ and $-11$ do not occur in this table. This is because $3$ is not inert in them and hence they do not give rise to an integral point on~$X_{ns}^{+}(9)$. There cannot be any imaginary quadratic field of class number one other than these nine fields, because in an imaginary quadratic field of class number one with discriminant $d$, all primes less than ${1+|d|}\over{4}$ are inert;    see $[{\bf 9}, {\rm{p}}.\,190]$. Moreover, the final entry in Table $5.2$ is a $j$-invariant that is not CM.  Hence if there were to be any other imaginary quadratic field with class number $1$, it would give rise to an integral point not on our exhaustive list.

\smallskip
We consider the last entry in Table $5.2$ in some more detail. It is the $j$-invariant $1117947^{{\hbox{$ \scriptstyle3$}}}
=3^{{\hbox{$ \scriptstyle3$}}}41^{{\hbox{$ \scriptstyle3$}}}61^{{\hbox{$ \scriptstyle3$}}}149^{{\hbox{$ \scriptstyle3$}}}$. A Weierstrass equation of an elliptic curve $E$ with this $j$-invariant is
$$
y^{2}+xy+y=x^{3}-x^{2}-408865825x-3182038133498.
$$
The discriminant of the equation of $E$ is equal to $5^{3}3511^{3}$. Counting points modulo the primes $p$ of good reduction with  $p<100$, one finds that $a_{p}(E)= p+1-|E(\FF_{p})|$ with $a_{p}$ as in Table $5.3$.

\bigskip
\vbox{\noindent {\bf Table 5.3.} 
\medskip
\centerline {\vbox {\offinterlineskip
\hrule\halign{&\vrule#&\strut\quid\hfil#\quid\cr
height3pt
&\omit&&\omit&&\omit&&\omit&&\omit&&\omit&&\omit&&\omit&&\omit&&\omit&&\omit&&\omit&&\omit&&\omit&&\omit&&\omit&&\omit&&\omit&&\omit&&\omit&&\omit&
&\omit&&\omit&&\omit&&\omit& \cr
&$p$&&2&&3&&7&&11&&13&&17&&19&&23&&29&&31&&37&&41&&43&&47&&53&&59&&61&&67&&71&&73&&79&&83&&89&&97&\cr
height3pt
&\omit&&\omit&&\omit&&\omit&&\omit&&\omit&&\omit&&\omit&&\omit&&\omit&&\omit&&\omit&&\omit&&\omit&&\omit&&\omit&&\omit&&\omit&&\omit&&\omit&&\omit&&\omit&&\omit&&\omit&&\omit& \cr
\noalign{\hrule}
height2pt
&\omit&&\omit&&\omit&&\omit&&\omit&&\omit&&\omit&&\omit&&\omit&&\omit&&\omit&&\omit&&\omit&&\omit&&\omit&&\omit&&\omit&&\omit&&\omit&&\omit&&\omit&&\omit&&\omit&&\omit&
&\omit&\cr
&$a_{p}^{}$&&$-1$&&0&&0&&$-2$&&0&&$-2$&&0&&$-1$&&0&&$-9$&&$-2$&&0&&1&&$-7$&&$-8$&&$-10$&&0&&$-4$&&9&&9&&8&&$-18$&&$-11$&&$-9$&\cr
height3pt
&\omit&&\omit&&\omit&&\omit&&\omit&&\omit&&\omit&&\omit&&\omit&&\omit&&\omit&&\omit&&\omit&&\omit&&\omit&&\omit&&\omit&&\omit&&\omit&&\omit&&\omit&&\omit&&\omit&&\omit&&\omit&\cr
}\hrule}}}
\bigskip\noindent
Using the techniques in $[{\bf{10}}]$ (especially Lemma $3$, p.$296$), we find that the quadratic subfield of $\QQ(E[9])$ that is the invariant field of the non-split Cartan subgroup of index $2$ is~$\QQ(\sqrt{-3511})$. The class number of $\QQ(\sqrt{-3511})$ is not $1$, but it is $41$. Half of the primes $p$, more precisely, the ones that are inert in $\QQ(\sqrt{-3511})$, satisfy $a_{p} \equiv 0 \zmod{9}$. This reflects the fact that the Galois group acts on $E[9]$ through the normalizer of a non-split Cartan subgroup. In some sense, this elliptic curve behaves as if it admits complex multiplication modulo $9$. 
\beginsection Acknowledgement

I would like to thank Ren\'e Schoof for introducing me to this problem and for his constant support and also for his patience in teaching me how to write a mathematics paper.

\bigskip\bigskip\noindent
$$\bf\rm{REFERENCES}$$

\item{[{\bf{1}}]} Baker, A., {\sl Linear forms in the logarithms of algebraic numbers}, Mathematika {\bf13} (1966), 204--216.
\smallskip
\item{[{\bf{2}}]} Baran, B., {\sl Normalizers of non-split Cartan subgroups, modular curves and the class number one problem}, in preparation.
\smallskip
\item{[{\bf{3}}]} Chen, I., 
{\sl On Siegel's modular curve of level $5$ and the class number one problem}, Journal of Number Theory {\bf 74} (1999), no.2, 278--297.
\smallskip
\item{[{\bf{4}}]} Deuring, M., {\sl Imagin\"are quadratische Zahlk\"orper mit der Klassenzahl Eins}, Inventiones Math. {\bf5} (1969), 169--179.
\smallskip
\item{[{\bf{5}}]} Diamond, F., Shurman, J., {\sl A First Course in Modular Forms}, Graduate Texts in Mathematics {\bf 228} (2005), Springer.
\smallskip
\item{[{\bf{6}}]} Heegner, K., {\sl Diophantische Analysis und Modulfunktionen}, Math. Zeit. {\bf59} (1952), 227--253.
\smallskip
\item{[{\bf{7}}]} Kenku, M.A., {\sl A note on the integral points of a modular curve of level $7$}, Mathematika {\bf 32} (1985), 45--48.
\smallskip
\item{[{\bf{8}}]} Ljunggren, W., {\sl Einige Bemerkungen \"uber die Darstellung ganzer Zahlen durch bin\"are kubische Formen mit positiver Diskriminante}, Acta Mathematica {\bf75} (1943), 1--21.
\smallskip
\item{[{\bf{9}}]} Serre, J.P., {\sl Lectures on the Mordell-Weil Theorem}, Aspects of Mathematics {\bf{E15}} (1989), Vieweg, Braunschweig/Wiesbaden.
\smallskip
\item{[{\bf{10}}]} Serre, J.P., {\sl Propri\'et\'es galoisiennes des points d'ordre fini des courbes elliptiques}, Inventiones Mathematicae {\bf{15}} (1972), 259--331.
\smallskip
 \item{[{\bf{11}}]} Shimura, G., {\sl Introduction to the Arithmetic Theory of Automorphic Functions}, Princeton University Press and Iwanami Shoten (1971), Princeton-Tokyo.
\smallskip
\item{[{\bf{12}}]} Siegel, C.L., {\sl Zum Beweise des Starkschen Satzes}, Inventiones Mathematicae {\bf{5}} (1968), 180--191.
\smallskip
\item{[{\bf{13}}]} Stark, H.M., {\sl On complex quadratic fields with class number equal to one}, Trans. Amer. Math. Soc. {\bf122} (1966), 112--119.
\smallskip
\item{[{\bf{14}}]} Stark, H.M., {\sl On the ``Gap" in a theorem of Heegner}, Journal of Number Theory {\bf1} (1969), 16--27.

\bye